\documentclass[11pt]{amsart}
\usepackage{amssymb}
\usepackage{amsmath}
\usepackage[active]{srcltx}
\usepackage{t1enc}
\usepackage[latin2]{inputenc}
\usepackage{verbatim}
\usepackage{amsmath,amsfonts,amssymb,amsthm}
\usepackage[mathcal]{eucal}
\usepackage{enumerate}
\usepackage[centertags]{amsmath}
\usepackage{graphics}

\newtheorem{theorem}{Theorem}

\newtheorem{corollary}{Corollary}

\begin{document}
\author{G. Tutberidze}
\title[partial sums]{A note on the strong convergence of partial
sums with respect to Vilenkin system}
\address{G.Tutberidze, The University of Georgia, school of Informatics, Engineering and
Mathematics, IV, 77a Merab Kostava St,
Tbilisi, 0128, Georgia, \& Department of Engineering Sciences and
Mathematics, Lule\aa University of Technology, SE-971 87 Lule\aa , Sweden.}
\email{giorgi.tutberidze1991@gmail.com}
\thanks{The research was supported by a Swedish Institute scholarship and by
Shota Rustaveli National Science Foundation grant YS15-2.1.1-47.}
\date{}
\maketitle

\begin{abstract}
In this paper we investigate some strong convergence theorems for partial
sums with respect to Vilenkin system.
\end{abstract}

\date{}

\textbf{2010 Mathematics Subject Classification.} 42C10.

\noindent \textbf{Key words and phrases:} Vilenkin system, partial sums, Fejé%
r means, Hardy space, strong convergence.

\section{Introduction}

It is well-known (for details see e.g. \cite{gol} and \cite{sws}) that the
Vilenkin system does not form a basis in the space $L_{1}\left( G_{m}\right).$
Moreover, there is a function in the Hardy space $H_{1}\left( G_{m}\right),$ (for details see \cite{ptw1, ptw2, tep5, tep6}) such that the partial sums of $f$ \ are not bounded in $L_{1}$-norm.
However, (see e.g. \cite{BNPT,tep8}) the subsequence $S_{M_{n}}$ of partial sums are bounded from the Hardy
space $H_{1}\left( G_{m}\right) $ to the Lebesgue space $L_{1}\left(
G_{m}\right):$
\begin{equation}  \label{1ccss}
\left\Vert S_{M_k}f\right\Vert _{H_1}\leq c\left\Vert f\right\Vert _{H_{1}}
\text{ \ \ \ } (k\in \mathbb{N}).
\end{equation}

Moreover, in Gát \cite{gat1} (see also Simon \cite{Si4,Si5}) it was proved the following strong convergence result for
all $f\in H_{1}:$%
\begin{equation*}
\underset{n\rightarrow \infty }{\lim }\frac{1}{\log n}\overset{n}{\underset{%
k=1}{\sum }}\frac{\left\Vert S_{k}f-f\right\Vert _{1}}{k}=0,
\end{equation*}%
where $S_{k}f$ denotes the $k$-th partial sum of the Vilenkin-Fourier series
of $f.$

It follows that there exists an absolute constant $c,$ such that%
\begin{equation}  \label{si}
\frac{1}{\log n}\overset{n}{\underset{k=1}{\sum }}\frac{\left\Vert
S_{k}f\right\Vert _{1}}{k}\leq c\left\Vert f\right\Vert _{H_{1}} \text{ \ }%
\left( n=2,3...\right)
\end{equation}%
and
\begin{equation*}
\underset{n\rightarrow \infty }{\lim }\frac{1}{\log n}\overset{n}{\underset{%
k=1}{\sum }}\frac{\left\Vert S_{k}f\right\Vert _{1}}{k}=\left\Vert
f\right\Vert _{H_{1}},
\end{equation*}%
for all $f\in H_{1}.$

Analogical result for the trigonometric system was proved by Smith \cite{sm}%
, and for the Walsh-Paley system by Simon \cite{Si3}.

If the partial sums of Vilenkin-Fourier series was bounded from $H_{1}$ to $%
L_{1} $ we also would have:
\begin{equation}  \label{tut}
\underset{n\in \mathbf{\mathbb{N}}_{+}}{\sup }\frac{1}{n}\underset{m=1}{%
\overset{n}{\sum }}\left\Vert S_{m}f\right\Vert _{1}\leq c\left\Vert
f\right\Vert_{H_1},
\end{equation}
but as it was presented above the boundednes of the partial sums does not hold
from $H_{1}$ to $L_{1} $, However, we have inequality (\ref{si}).

On the other hand, in the one-dimensional case, Fujji \cite{Fu} and Simon \cite{Si2}
proved that maximal operator Fejér means are bounded from $H_{1}$ to $L_{1}$,
that is
\begin{equation}  \label{fusi}
\underset{n\in \mathbf{\mathbb{N}}_{+}}{\sup }\left\Vert \frac{1}{n}\underset%
{m=1}{\overset{n}{\sum }} S_{m}f\right\Vert _{1}<c\left\Vert f
\right\Vert_{H_1}.
\end{equation}

So, natural question has arised that if inequality (\ref{tut}) holds true,
which would be generalization of inequality (\ref{fusi}) or is we have negative
answer on this problem.

In this paper we prove that there exists a function $f\in H_{1} $ such that
\begin{equation*}
\underset{n\in \mathbf{\mathbb{N}}_{+}}{\sup }\frac{1}{n}\underset{m=1}{%
\overset{n}{\sum }}\left\Vert S_{m}f\right\Vert _{1}=\infty.
\end{equation*}

This paper is organized as follows: in order not to disturb our discussions
later on some definitions and notations are presented in Section 2. For the
proofs of the main results we need some auxiliary Lemmas. These results are
presented in Section 3. The formulation and detailed proof of main results
can be found in Section 4.

\section{Definitions and Notations}

Let $\mathbb{N}_{+}$ denote the set of the positive integers, $\mathbb{N}:=%
\mathbb{N}_{+}\cup \{0\}.$

Let $m:=(m_{0,}m_{1},\dots)$ denote a sequence of  positive integers not
less than 2.

Denote by
\begin{equation*}
Z_{m_{k}}:=\{0,1,\dots,m_{k}-1\}
\end{equation*}
the additive group of integers modulo $m_{k}.$

Define the group $G_{m}$ as the complete direct product of the group $%
Z_{m_{j}}$ with the product of the discrete topologies of $Z_{m_{j}}$ $^{,}$%
s.

The direct product $\mu $ of the measures
\begin{equation*}
\mu _{k}\left( \{j\}\right):=1/m_{k}\text{ \qquad }(j\in Z_{m_{k}})
\end{equation*}
is the Haar measure on $G_{m_{\text{ }}}$with $\mu \left( G_{m}\right) =1.$

If $\sup_{n\in \mathbb{N}}m_{n}<\infty $, then we call $G_{m}$ a bounded
Vilenkin group. If the generating sequence $m$ is not bounded, then $G_{m}$
is said to be an unbounded Vilenkin group.

The elements of $G_{m}$ are represented by sequences
\begin{equation*}
x:=(x_{0},x_{1},\dots,x_{k},\dots)\qquad \left( \text{ }x_{k}\in
Z_{m_{k}}\right).
\end{equation*}

It is easy to give a base for the neighbourhood of $G_{m}$ namely
\begin{eqnarray*}
I_{0}\left( x\right)&:=&G_{m}, \\
I_{n}(x)&:=&\{y\in G_{m}\mid y_{0}=x_{0},\dots,y_{n-1}=x_{n-1}\}\text{ }%
(x\in G_{m},\text{ }n\in \mathbb{N})
\end{eqnarray*}%
Denote $I_{n}:=I_{n}\left( 0\right) $ for $n\in \mathbb{N}$ and $\overline{%
I_{n}}:=G_{m}$ $\backslash $ $I_{n}$ $.$

Let
\begin{equation*}
e_{n}:=\left( 0,\dots,0,x_{n}=1,0,\dots\right) \in G_{m}\qquad \left( n\in
\mathbb{N}\right).
\end{equation*}

If we define the so-called generalized number system based on $m$ in the
following way:
\begin{equation*}
M_{0}:=1,\text{\qquad }M_{k+1}:=m_{k}M_{k\text{ }},\ \qquad (k\in \mathbb{N})
\end{equation*}%
then every $n\in \mathbb{N}$ can be uniquely expressed as $%
n=\sum_{k=0}^{\infty }n_{j}M_{j}$ where $n_{j}\in Z_{m_{j}}$ $~(j\in \mathbb{%
N})$ and only a finite number of $n_{j}`$s differ from zero. Let $\left\vert
n\right\vert :=\max $ $\{j\in \mathbb{N};$ $n_{j}\neq 0\}.$

Next, we introduce on $G_{m}$ an orthonormal system, which is called the
Vilenkin system.

At first define the complex valued function $r_{k}\left( x\right)
:G_{m}\rightarrow \mathbb{C},$ the generalized Rademacher functions as
\begin{equation*}
r_{k}\left( x\right):=\exp \left( 2\pi\imath x_{k}/m_{k}\right) \text{
\qquad }\left( \imath^{2}=-1,\text{ }x\in G_{m},\text{ }k\in \mathbb{N}%
\right).
\end{equation*}

Now define the Vilenkin system $\psi:=(\psi _{n}:n\in \mathbb{N})$ on $G_{m}
$ as:
\begin{equation*}
\psi _{n}\left( x\right):=\prod_{k=0}^{\infty }r_{k}^{n_{k}}\left( x\right)
\text{ \qquad }\left( n\in \mathbb{N}\right).
\end{equation*}

Specially, we call this system the Walsh-Paley one if $m\equiv 2.$

The norm (or quasi norm) of the space $L_{p}(G_{m})$ is defined by \qquad
\qquad \thinspace\
\begin{equation*}
\left\Vert f\right\Vert _{p}:=\left( \int_{G_{m}}\left\vert f(x)\right\vert
^{p}d\mu (x)\right) ^{1/p}\qquad \left( 0<p<\infty \right) .
\end{equation*}

The Vilenkin system is orthonormal and complete in $L_{2}\left( G_{m}\right)
\,$ (for details see e.g. \cite{AVD,Vi}).

If $\ f\in L_{1}\left( G_{m}\right) $ we can define Fourier coefficients,
partial sums of the Fourier series, Fejér means, Dirichlet kernels with
respect to the Vilenkin system in the usual manner:
\begin{eqnarray*}
\widehat{f}(k) &:&=\int_{G_{m}}f\overline{\psi }_{k}d\mu \text{\thinspace
\qquad\ \ \ \ }\left( \text{ }k\in \mathbb{N}\text{ }\right) \\
S_{n}f &:&=\sum_{k=0}^{n-1}\widehat{f}\left( k\right) \psi _{k}\ \text{%
\qquad\ \ }\left( \text{ }n\in \mathbb{N}_{+},\text{ }S_{0}f:=0\right) \\
\sigma _{n}f &:&=\frac{1}{n}\sum_{k=0}^{n-1}S_{k}f\text{ \qquad\ \ \ \ \ }%
\left( \text{ }n\in \mathbb{N}_{+}\text{ }\right) \\
D_{n} &:&=\sum_{k=0}^{n-1}\psi _{k\text{ }}\text{ \qquad\ \ \qquad }\left(
\text{ }n\in \mathbb{N}_{+}\text{ }\right).
\end{eqnarray*}

Recall that
\begin{equation}  \label{3}
\quad \hspace*{0in}D_{M_{n}}\left( x\right) =\left\{
\begin{array}{l}
\text{ }M_{n},\text{ \ \ \ }x\in I_{n} \\
\text{ }0,\text{ \qquad }x\notin I_{n}%
\end{array}
\right.
\end{equation}
and
\begin{equation}  \label{9dn}
D_{s_{n}M_{n}}=D_{M_{n}}\sum_{k=0}^{s_{n}-1}\psi
_{kM_{n}}=D_{M_{n}}\sum_{k=0}^{s_{n}-1}r_{n}^{k}, \text{ \qquad } 1\leq
s_{n}\leq m_{n}-1.
\end{equation}

The $n$-th Lebesgue constant is defined in the following way
\begin{equation*}
L_{n}=\left\Vert D_{n}\right\Vert _{1}.
\end{equation*}

It is well-known \cite{Vi} that 
\begin{equation} \label{dnl}
L_{n}=O(\log n), \ \ n\rightarrow \infty.
\end{equation}

Moreover, (for unbounded Vilenkin systems it can be found in \cite{FS},
for bounded Vilenkin systems see e.g. \cite{luko} and \cite{smt,tep9}) there exist absolute constant 
$c_1 $ and $c_2 $ such that
\begin{equation} \label{estln}
c_1\log n\leq\frac{1}{n}\underset{k=1}{\overset{n}{\sum }}L\left( k\right)
\leq c_2 \log n, \ \ (n=2,3,...).
\end{equation}

The concept of the Hardy space \cite{CW} can be defined in various manners, e.g.
by a maximal function 
\[
f^{\ast }:=\sup_{n}\left\vert S_{M_{n}}f\right\vert \text{ \ \ }\left( f\in
G_{m}\right) ,
\]%
saying that $f$ belongs to the Hardy space if  $f^{\ast }\in L^{1}\left(
G_{m}\right) .$ This definition is suitable if the sequence $m$ is bounded.
In this case a good property of the space $\left\{ f\in L^{1}\left(
G_{m}\right) :\text{ }f^{\ast }\in L^{1}\left( G_{m}\right) \right\} $ is
the atomic structure \cite{CW}. To the definition of space of Hardy type for
an arbitrary $m$, first we give the concept of the atoms \cite{Si2}. A set $%
I\subset G_{m}$ is called an interval if for some $x\in G_{m}$ and $n\in N$, 
$I$ is of the form  $I=\bigcup\limits_{k\in U}I_{n}\left(x,k\right)$, where  $U$ is obtained from $Z_{m_{n}}$ by dyadic partition.

The sets $U_{1},U_{2},...\subset Z_{m_{n}}$, are obtained by means of such
a partition if%
\[
U_{1}=\left\{ 0,...,\left[ \frac{m_{n}}{2}\right] -1\right\} ,\text{ \ }%
U_{2}=\left\{ \left[ \frac{m_{n}}{2}\right] ,...,m_{n}-1\right\} ,
\]%
\[
U_{3}=\left\{ 0,...,\left[ \frac{\left[ m_{n}/2\right] -1}{2}\right]
-1\right\} ,\text{ \ }U_{4}=\left\{ \left[ \frac{[m_{n}/2]-1}{2}\right] ,...,%
\left[ \frac{m_{n}}{2}\right] -1\right\} ,...
\]%
etc.; $[\ \ ]$ denotes the entire part. We define the atoms as follows: the
function $a\in L^{\infty }\left( G_{m}\right) $ is called an atom if
eather $a\equiv 1$ or there exists an interval $I$ for which $\sup a\subset I,$ $%
\left\vert a\right\vert \leq \left\vert I\right\vert ^{-1}$ and $\int_{I}a=0$
hold. $\left( \left\vert I\right\vert \text{ denotes the Haar measure of }I%
\text{ }\right) .$ 

Now we can define the space $H_1\left( G_{m}\right) $ (for details see e.g \cite{We1, We5}) as the set of all functions $f=\sum\limits_{i=0}^{\infty }\lambda _{i}a_{i}$,where $a_{i}$'s are atoms and for the coefficients $\lambda _{i}$ we have $%
\sum\limits_{i=0}^{\infty }\left\vert \lambda _{i}\right\vert <\infty .$ $%
H_1\left( G_{m}\right) $ is a Banach space with respect to the norm 
\begin{equation} \label{norm}
\left\Vert f\right\Vert_{H_1}:=\inf \sum\limits_{k=0}^{\infty }\left\vert
\lambda _{k}\right\vert <\infty.
\end{equation}

The infimum is taken over all decompositions 
$$ f=\sum\limits_{i=0}^{\infty }\lambda _{i}a_{i}.$$ It is known [7]
that $\left\Vert f\right\Vert_{H_1} $ is equivalent to $\left\Vert f^{\ast \ast
}\right\Vert _{1}\left( f\in L^{1}\left( G_{m}\right) \right) $, where $f$**$%
\left( x\right) :=\sup_{I}\left\vert I\right\vert ^{-1}\left\vert
\int_{I}f\right\vert $, $\left( x\in G_{m}\text{, }x\in I\text{ and }I\text{
is interval}\right) $. Since by (\ref{3})
\begin{equation*}
f^{*}\left( x\right) =\sup_{n\in \mathbb{N}}\frac{1}{\left| I_{n}\left(
x\right) \right| }\left| \int_{I_{n}\left( x\right) }f\left( u\right) \mu
\left( u\right) \right|
\end{equation*}
we have $f^{\ast }\leq f^{\ast \ast }$ and, thus, $H\left( G_{m}\right)
\subset \left\{ f\in L^{1}\left( G_{m}\right) :\text{ }f^{\ast }\in
L^{1}\left( G_{m}\right) \right\} .$ Moreover these spaces coincide if the sequence $m$ is bounded.

\section{Main Result}
Our main result reads:
\begin{theorem}
\label{theorem1}a) Let $f\in H_{1}.$ Then there exists an absolute constant $%
c,$ such that
\begin{equation*}
\sup_{n\in \mathbb{N}}\frac{1}{n\log n}\overset{n}{\underset{k=1}{\sum }}%
\left\Vert S_{k}f\right\Vert _{1}\leq \left\Vert f\right\Vert_{H_1}.
\end{equation*}

b) Let $\varphi :N_{+}\rightarrow \lbrack 1,$ $\infty )$ be a nondecreasing
function satisfying the condition
\begin{equation}
\overline{\lim_{n\rightarrow \infty }}\frac{\log n}{\varphi _{n}}=+\infty .
\label{cond1}
\end{equation}

Then there exists a function $f\in H_{1},$ such that
\begin{equation*}
\sup_{n\in \mathbb{N}}\frac{1}{n\varphi_n}\overset{n}{\underset{k=1}{\sum }}%
\left\Vert S_{k}f\right\Vert _{1}=\infty .
\end{equation*}
\end{theorem}

\begin{corollary}
\label{theorem2} (see e.g. \cite{MS,Si2,Si4})There exists a function $f\in H_{1},$ such that
\begin{equation*}
\sup_{n\in\mathbb{N}}\frac{1}{n}\overset{n}{\underset{k=1}{\sum }}\left\Vert
S_{k}f\right\Vert_{1}=\infty .
\end{equation*}
\end{corollary}

\section{Proof of Theorem \protect\ref{theorem2}}

\begin{proof}[Proof]
By using (\ref{dnl}) we can conclude that 
\begin{equation*}
\frac{1}{n\log n}\overset{n}{\underset{k=1}{\sum }}%
\left\Vert S_{k}f\right\Vert _{1}\leq  \frac{c\left\Vert f\right\Vert_{H_1}}{n\log n}\overset{n}{\underset{k=1}{\sum }}{\log k}\leq c\left\Vert f\right\Vert_{H_1}.
\end{equation*}

The proof of part a) is complete.

Under the condition \eqref{cond1} there exists an increasing sequence of the positive integers $\left\{\alpha_{k}:k\in \mathbb{N}\right\} $ such that
\begin{equation*}
\overline{\lim_{k\rightarrow \infty }}\frac{\log M_{{\alpha _{k}}}}{\varphi
_{2M_{\alpha _{k}}}}=+\infty
\end{equation*}%
and
\begin{equation} \label{69}
\sum_{k=0}^{\infty }\frac{\varphi _{2M_{\alpha _{k}}}^{1/2}}{\log ^{1/2}M_{{%
\alpha _{k}}}}<c<\infty .  
\end{equation}

Let \qquad
\begin{equation*}
f=\sum_{k=1}^{\infty }{\lambda _{k}}{a_{k}},
\end{equation*}%
where
\begin{equation*}
a_{k}=r_{\alpha _{k}}D_{M_{_{\alpha _{k}}}}=D_{2M_{_{\alpha
_{k}}}}-D_{M_{_{\alpha _{k}}}}
\end{equation*}%
and
\begin{equation*}
\lambda _{k}=\frac{\varphi _{2M_{\alpha _{k}}}^{1/2}}{\log ^{1/2}M_{{\alpha
_{k}}}}.
\end{equation*}

By the definition of $H_1$ and (\ref{norm}), if we apply (\ref{69}) we can conclude that $ f\in H_{1}.$ Moreover,
\begin{equation} \label{6}
\widehat{f}(j)=\left\{
\begin{array}{l}
{\lambda _{k}},\,\,\text{\ \ \ \ \thinspace \thinspace }j\in \left\{
M_{\alpha _{k}},...,2M_{\alpha _{k}}-1\right\} ,\text{ }k\in \mathbb{N} \\
0\text{ },\text{ \thinspace \qquad \thinspace \thinspace \thinspace
\thinspace \thinspace }j\notin \bigcup\limits_{k=1}^{\infty }\left\{
M_{\alpha _{k}},...,2M_{\alpha _{k}}-1\right\} .\text{ }%
\end{array}%
\right.   
\end{equation}

Since
\begin{equation*}
D_{j+M_{\alpha _{k}}}=D_{M_{\alpha _{k}}}+\psi _{_{M_{\alpha _{k}}}}D_{j},%
\text{ \qquad when \thinspace \thinspace }j\leq M_{\alpha _{k}},
\end{equation*}
if we apply (\ref{6}) we obtain that
\begin{eqnarray}
S_{j}f &=&S_{M_{\alpha _{k}}}f+\sum_{v=M_{\alpha _{k}}}^{j-1}\widehat{f}%
(v)\psi _{v}  \label{8} \\
&=&S_{M_{\alpha _{k}}}f+{\lambda _{k}}\sum_{v=M_{\alpha _{k}}}^{j-1}\psi _{v}
\notag \\
&=&S_{M_{\alpha _{k}}}f+{\lambda _{k}}\left( D_{j}-D_{M_{\alpha _{k}}}\right)
\notag \\
&=&S_{M_{\alpha _{k}}}f+{\lambda _{k}}\psi _{M_{\alpha _{k}}}D_{j-M_{\alpha
_{k}}}  \notag \\
&=&I_{1}+I_{2}.  \notag
\end{eqnarray}

In view of (\ref{1ccss}) we can write that
\begin{eqnarray}  \label{8bbb2}
\left\Vert I_{1}\right\Vert_{1} \leq\left\Vert
S_{M_{\alpha_{k}}}f\right\Vert _{1} \leq c\left\Vert f\right\Vert_{H_{1}}.
\end{eqnarray}

By combining (\ref{estln}) and (\ref{8bbb2}) we get that%
\begin{eqnarray*}
\left\Vert S_{n}f\right\Vert _{1} \geq \left\Vert I_{2}\right\Vert
_{1}-\left\Vert I_{1}\right\Vert _{1} \geq {\lambda _{k}}{L\left( {%
n-M_{\alpha _{k}}}\right)}-c\left\Vert f\right\Vert _{H_{1}}.
\end{eqnarray*}

Hence, 
\begin{eqnarray*}
&&\underset{n\in \mathbf{\mathbb{N}}_{+}}{\sup }\frac{1}{n\varphi _{n}}%
\underset{k=1}{\overset{n}{\sum }}\left\Vert S_{k}f\right\Vert _{1} \\
&\geq &\frac{1}{2M_{\alpha _{k}}\varphi _{2M_{\alpha _{k}}}}\underset{%
\left\{ M_{\alpha _{k}}\leq l\leq 2M_{\alpha _{k}}\right\} }{\sum }%
\left\Vert S_{l}f\right\Vert _{1} \\
&\geq &\frac{1}{2M_{\alpha _{k}}\varphi _{2M_{\alpha _{k}}}}\underset{%
\left\{ M_{\alpha _{k}}\leq l\leq 2M_{\alpha _{k}}\right\} }{\sum }\left(
\frac{L\left( l-M_{\alpha _{k}}\right) \varphi _{2M{\alpha _{k}}}^{1/2}}{%
\log ^{1/2}M_{\alpha _{k}}}-c\left\Vert f\right\Vert _{H_{1}}\right) \\
&\geq &\frac{c\varphi _{2M{\alpha _{k}}}^{1/2}}{2M_{\alpha _{k}}\log
^{1/2}M_{\alpha _{k}}\varphi _{2M_{\alpha _{k}}}}\underset{l=1}{\overset{%
M_{\alpha _{k}}-1}{\sum }}L\left( l\right) -c\left\Vert f\right\Vert
_{H_{1}}^{1/2} \\
&\geq &\frac{c\varphi _{2M{\alpha _{k}}}^{1/2}\log M_{\alpha _{k}}}{\log
^{1/2}M_{\alpha _{k}}\varphi _{2M_{\alpha _{k}}}}\geq \frac{c\log
^{1/2}M_{\alpha _{k}}}{\varphi _{2M_{\alpha _{k}}}^{1/2}}\rightarrow \infty ,%
\text{ as \ }k\rightarrow \infty .
\end{eqnarray*}

The proof is complete.
\end{proof}

\textbf{Acknowledgment:} The author would like to thank the referee for helpful suggestions, which improved the final version of the paper.

\end{document}